\documentclass[a4paper,11pt]{article}
\usepackage{amssymb}

\newcommand{\R}{\mathbb{R}}

\newcommand{\beq}{\begin{equation} }
\newcommand{\eqq}{\end{equation} }
\newcommand{\cuad}{{\sqcap\kern-.68em\sqcup}}
\newcommand{\abs}[1]{\mid #1 \mid}
\newcommand{\norm}[1]{\|#1\|}

\newtheorem{definition}{Definition}[section]
\newtheorem{teo}{Theorem}[section]

\newtheorem{proposition}{Proposition}[section]

\newtheorem{lemma}{Lemma}[section]
\newtheorem{corollary}{Corollary}[section]
\newtheorem{remark}{Remark}[section]
\newcommand{\bremark}{\begin{remark} \em}
\newcommand{\eremark}{\end{remark} }

\def\beeq{\begin{equation}}
\def\eeq{\end{equation}}
\newcommand{\begeqaet}{\begin{eqnarray*}}
\newcommand{\eneqaet}{\end{eqnarray*}}

\newcommand{\re}[1]{(\ref{#1})}
\begin{document}

\begin{center}{\bf  \Large   Weakly and strongly singular solutions of\\[2mm] semilinear fractional  elliptic equations }\medskip
\medskip

\bigskip
{\small
      {\bf Huyuan Chen}\footnote{chenhuyuan@yeah.net}\smallskip
      
       Departamento de Ingenier\'{\i}a  Matem\'atica
\\ Universidad de Chile, Chile\\[2mm]
     {\bf Laurent V\'{e}ron}\footnote{Laurent.Veron@lmpt.univ-tours.fr}

\smallskip
Laboratoire de Math\'{e}matiques et Physique Th\'{e}orique
\\  Universit\'{e} Fran\c{c}ois Rabelais, Tours, France}\\[1mm]
\bigskip

\medskip

\begin{abstract}
Let $p\in(0,\frac{N}{N-2\alpha})$, $\alpha\in(0,1)$ and $\Omega\subset \R^N$ be a bounded $C^2$ domain containing $0$. If $\delta_0$ is the Dirac measure at $0$ and $k>0$, we prove that  the  weakly singular solution $u_k$ of $(E_k)$ $ (-\Delta)^\alpha
u+u^p=k\delta_0 $ in $\Omega$ which vanishes in $\Omega^c$, is a classical  solution of
$(E_*)$ $ (-\Delta)^\alpha
u+u^p=0 $ in $\Omega\setminus\{0\}$ with the same outer data. When $\frac{2\alpha}{N-2\alpha}\leq 1+\frac{2\alpha}{N}$, $p\in(0, 1+\frac{2\alpha}{N}]$ we show that the $u_k$ converges to $\infty$ in whole $\Omega$ when $k\to\infty$, while, for $p\in(1+\frac{2\alpha}N,\frac{N}{N-2\alpha})$, the limit of the $u_k$ is a strongly singular solution of $(E_*)$. The same result holds in the case $1+\frac{2\alpha}{N}<\frac{2\alpha}{N-2\alpha}$ excepted if $\frac{2\alpha}{N}<p<1+\frac{2\alpha}{N}$.
\end{abstract}
\end{center}
\tableofcontents \vspace{1mm}
  \vspace{1mm}


\vspace{2mm}
\setcounter{equation}{0}
\section{Introduction}

Let $\Omega$ be a bounded  $C^2$ domain of $\R^N(N\geq2)$ containing $0$, $ \alpha\in(0,1)$ and let $\delta_0$ denote the Dirac measure at $0$.
In this paper, we study the properties of the weak solution to problem
\begin{equation}\label{eq1.1}
 \arraycolsep=1pt
\begin{array}{lll}
 (-\Delta)^\alpha  u+u^p=k\delta_0\quad & \rm{in}\quad\Omega\\[2mm]
 \phantom{   (-\Delta)^\alpha  + u^p}
u=0\quad & \rm{in}\quad \Omega^c,
\end{array}
\end{equation}
where  $k>0$ and $p\in(0,\frac N{N-2\alpha})$ and $(-\Delta)^\alpha $
  is the $\alpha$-fractional Laplacian defined by
$$(-\Delta)^\alpha  u(x)=\lim_{\epsilon\to0^+} (-\Delta)_\epsilon^\alpha u(x),$$
where for $\epsilon>0$,
$$
(-\Delta)_\epsilon^\alpha  u(x)=-\int_{\R^N}\frac{ u(z)-
u(x)}{|z-x|^{N+2\alpha}}\chi_\epsilon(|x-z|) dz
$$
and
$$\chi_\epsilon(t)=\left\{ \arraycolsep=1pt
\begin{array}{lll}
0,\quad & \rm{if}\quad t\in[0,\epsilon]\\[2mm]
1,\quad & \rm{if}\quad t>\epsilon.
\end{array}
\right.$$

In 1980, Benilan and Brezis (see \cite{Br, BB}) studied the case $\alpha=1$ in equation (\ref{eq1.1}) and proved in particular that  equation
\begin{equation}\label{eq003}
 \arraycolsep=1pt
\begin{array}{lll}
 -\Delta  u+u^q=k\delta_0 \quad & \rm{in}\quad\Omega\\[2mm]
 \phantom{ -\Delta  +u^q}
u=0  \quad & \rm{on}\quad \partial\Omega
\end{array}
\end{equation}
admits a unique solution $u_k$ for $1< q< N/(N-2)$, while no solution exists when $q\geq  N/(N-2)$. Soon after,
 Brezis and V\'{e}ron \cite{BV} proved that the problem
 \begin{equation}\label{local}
 \arraycolsep=1pt
\begin{array}{lll}
 -\Delta  u+u^q=0 \quad & \rm{in}\quad\Omega\setminus\{0\}\\[2mm]
 \phantom{ -\Delta  +u^q }
u=0  \quad & \rm{on}\quad \partial\Omega
\end{array}
\end{equation}
admits only the zero solution when $q\geq  N/(N-2)$. When $1< q< N/(N-2)$, V\'{e}ron in \cite{V0} obtained the description of the  all the possible singular behaviour of the positive  solutions of (\ref{local}). In particular he proved that this behaviour is always isotropic (when
$(N+1)/(N-1)\leq q< N/(N-2)$ the assumption of positivity is unnecessary) and that two types of singular behaviour occur:

\noindent (i) either $u(x)\sim c_Nk|x|^{2-N}$ when $x\to 0$ and $k$ can take any positive value; $u$ is said to have a {\it weak singularity} at $0$, and actually $u=u_k$.

\noindent (ii) or $u(x)\sim c_{N,q}|x|^{-\frac2{q-1}}$ when $x\to 0$ and $u$ has a {\it strong singularity} at $0$, and
$u=u_\infty:=\lim_{k\to\infty}u_k$.\smallskip

A large series of papers has been devoted to the extension of semilinear problems involving the Laplacian to problems where the diffusion operator is non-local, the most classical one being  the fractional Laplacian, see e.g. \cite{CS1, CFQ, CS, RS, S}. 
In a recent work, Chen and V\'{e}ron \cite{CV1} considered the
 problem
\begin{equation}\label{eq1.2}
 \arraycolsep=1pt
\begin{array}{lll}
 (-\Delta)^\alpha  u+u^p=0\quad & \rm{in}\quad\Omega\setminus\{0\}\\[2mm]
 \phantom{   (-\Delta)^\alpha  + u^p}
u=0\quad & \rm{in}\quad \Omega^c,
\end{array}
\end{equation}
where $1+\frac{2\alpha}N<p<p^*_{\alpha}:=\frac{N}{N-2\alpha}$.
They proved that (\ref{eq1.2}) admits a singular solution $u_s$ which satisfies
\begin{equation}\label{eq1.01}
\lim_{x\to0}u_s(x)|x|^{\frac{2\alpha}{p-1}}=c_0,
\end{equation}
for some $c_0>0$.
Moreover $u_s$  is the unique positive solution of (\ref{eq1.2}) such that
\begin{equation}\label{strong singular unique}
0<\liminf_{x\to0}u(x)|x|^{\frac{2\alpha}{p-1}}\le \limsup_{x\to0}u(x)|x|^{\frac{2\alpha}{p-1}}<\infty.
\end{equation}

In this article we will call {\it weakly singular solution} a solution $u$ of (\ref{eq1.2}) which satisfies $\limsup_{x\to0}|u(x)||x|^{N-2\alpha}<\infty$  and {\it strongly singular solution}
if  $\lim_{x\to0}|u(x)||x|^{N-2\alpha}=\infty.$

The existence of solutions of (\ref{eq1.1}) is a particular case of the more general problem  
\begin{equation}\label{eq1.01hc}
\arraycolsep=1pt
\begin{array}{lll}
 (-\Delta)^\alpha  u+g(u)=\nu\quad & \rm{in}\quad\Omega\\[2mm]
 \phantom{  (-\Delta)^\alpha  + g(u)}
u=0\quad & \rm{in}\quad \Omega^c
\end{array}
\end{equation}
which has been study by Chen and V\'{e}ron in \cite{CV2}  under the assumption that $g$ is a subcritical nonlinearity, $\nu$ being a positive and bounded Radon measure in $\Omega$. 

\begin{definition}\label{weak definition}
A function $u$ belonging to
$L^1(\Omega)$ is a weak solution of (\ref{eq1.01hc}) if $g(u)\in L^1(\Omega,\rho^\alpha dx)$   and
\begin{equation}\label{weak sense}
\int_\Omega [u(-\Delta)^\alpha\xi+g(u)\xi]dx=\int_\Omega \xi d\nu\qquad \forall\xi\in \mathbb{X}_{\alpha},
\end{equation}
where $\rho(x):=dist (x,\Omega^c)$ and $\mathbb{X}_{\alpha}\subset C(\R^N)$ is the space of functions
$\xi$ satisfying:\smallskip

\noindent (i) $\rm{supp}(\xi)\subset\bar\Omega$,\smallskip

\noindent(ii) $(-\Delta)^\alpha\xi(x)$ exists for all $x\in \Omega$
and $|(-\Delta)^\alpha\xi(x)|\leq c_1$ for some $c_1>0$,\smallskip

\noindent(iii) there exist $\varphi\in L^1(\Omega,\rho^\alpha dx)$
and $\epsilon_0>0$ such that $|(-\Delta)_\epsilon^\alpha\xi|\le
\varphi$ a.e. in $\Omega$, for all
$\epsilon\in(0,\epsilon_0]$.\smallskip
\end{definition}

According to Theorem 1.1 in \cite{CV2},  problem (\ref{eq1.1})
admits a unique weak solution $u_k$, moreover,
\begin{equation}\label{1.3}
  \mathbb{G}_\alpha[k\delta_0]-\mathbb{G}_\alpha[(\mathbb{G}_\alpha[k\delta_0])^p]\le u_k\le   \mathbb{G}_\alpha[k\delta_0]\quad {\rm in}\ \ \Omega,
\end{equation}
where $\mathbb{G}_\alpha[\cdot]$ is the Green operator
defined by
\begin{equation}\label{optimal0}
\mathbb{G}_\alpha[\nu](x)=\int_{\Omega}G_\alpha(x,y) d\nu(y),\qquad\forall\ \nu\in
\mathfrak{M}(\Omega,\rho^\alpha),
\end{equation}
with $G_\alpha$ is the Green kernel of $(-\Delta)^\alpha$ in
$\Omega\time\Omega$ and $\mathfrak{M}(\Omega,\rho^\alpha)$ denotes the
space of Radon measures in $\Omega$ such that $\int_{\Omega}\rho^\alpha d\abs\nu<\infty$.
By (\ref{1.3}),
\begin{equation}\label{1.4}
\lim_{x\to0}u_k(x)|x|^{N-2\alpha}=c_{\alpha,N}k.
\end{equation}
for some $c_{\alpha,N}>0$. From Theorem 1.1 in  \cite{CV2}, there holds
\begin{equation}\label{monotonicity}
 u_{k}(x)\le u_{k+1}(x),\qquad \forall x\in \Omega;
\end{equation}
then there exists
\begin{equation}\label{definition infty}
  u_{\infty}(x)=\lim_{k\to\infty}u_k(x) \quad\forall  x\in\R^N\setminus\{0\},
\end{equation}
and $u_{\infty}(x)\in \R_+\cup\{+\infty\}$.

Motivated by these results and in view of the nonlocal character of the fractional Laplacian, in this article we analyse  the connection between the solutions of (\ref{eq1.1}) and
the ones of (\ref{eq1.2}). Our main result is the following
\begin{teo}\label{teo 1}
Assume that $1+\frac{2\alpha}N\ge\frac{2\alpha}{N-2\alpha}$ and
$p\in(0,p^*_{\alpha})$.
Then
$u_k$ is a classical solution of (\ref{eq1.2}).
Furthermore, \\
  $(i)$ if $p\in(0,1+\frac{2\alpha}N)$,
\begin{equation}\label{xuan 10}
 u_\infty(x)=\infty \qquad \forall x\in \Omega;
\end{equation}
$(ii)$ if $p\in(1+\frac{2\alpha}N, p^*_\alpha)$,
 $$
   u_\infty =u_s,
$$
where $u_s$ is the solution of (\ref{eq1.2}) satisfying (\ref{eq1.01}).

\noindent Moreover, if $1+\frac{2\alpha}N=\frac{2\alpha}{N-2\alpha}$, (\ref{xuan 10}) holds
for  $p=1+\frac{2\alpha}N$.

\end{teo}

The result of part $(i)$ indicates that {\it even if the absorption is superlinear}, the diffusion dominates and there is no strongly singular solution to problem (\ref{eq1.2}).
On the contrary, part $(ii)$ points out that the absorption dominates the diffusion; the limit function $u_s$ is the  least  strongly singular solution of (\ref{eq1.2}).  
Comparing Theorem \ref{teo 1} with the results for Laplacian case, part $(i)$ with $p\in(0,1]$ and $(ii)$ are similar as the Laplacian case,
but part $(i)$ with $p\in(1,1+\frac{2\alpha}N]$ is totally different from the one in the case  $\alpha=1$. 
This striking phenomenon comes comes from the fact that the fractional Laplacian is a nonlocal operator, which requires the solution to belong to $L^1(\Omega)$, therefore no local barrier can be constructed if $p$ is too close to 1.\smallskip

At end, we consider the case where $1+\frac{2\alpha}N< \frac{2\alpha}{N-2\alpha}$. It occurs when
 $N=2$ and $ \frac{\sqrt{5}-1}2<\alpha<1$  or  $N=3$ and $\frac{3(\sqrt{5}-1)}4<\alpha<1$.
In this situation, we have the following  results.

\begin{teo}\label{teo 1.1}
Assume that $1+\frac{2\alpha}N< \frac{2\alpha}{N-2\alpha}$ and  $p\in(0, p^*_\alpha)$.
Then
$u_k$ is a classical solution of (\ref{eq1.2}).
Furthermore, \\
  $(i)$ if $p\in(0,\frac N{2\alpha})$, then
$$
    u_\infty(x)=\infty \qquad \forall x\in \Omega;
$$
$(ii)$ if $p\in( 1+\frac{2\alpha}N, \frac{2\alpha}{N-2\alpha})$, then
$u_\infty$ is a classical solution of (\ref{eq1.2}) and there exist $\rho_0>0$ and $c_2>0$ such that
\begin{equation}\label{1.2.1}
 c_2|x|^{-\frac{(N-2\alpha)p}{p-1}}\le u_\infty\le u_s \qquad \forall x\in B_{\rho_0}\setminus\{0\};
 \end{equation}
$(iii)$ if $p=\frac{2\alpha}{N-2\alpha}$, then
$u_\infty$ is a classical solution of (\ref{eq1.2}) and
there exist $\rho_0>0$ and $c_3>0$ such that
\begin{equation}\label{1.1.1}
  c_3\frac{|x|^{-\frac{(N-2\alpha)p}{p-1}}}{(1+|\log(|x|)|)^{\frac1{p-1}}}\le u_\infty\le u_s \qquad \forall x\in B_{\rho_0}\setminus\{0\};
 \end{equation}
$(iv)$ if $p\in(\frac{2\alpha}{N-2\alpha},p^*_\alpha)$, then
$$
u_\infty= u_s.
 $$
\end{teo}

We remark that $\frac N{2\alpha}<1+\frac{2\alpha}N$ if  $1+\frac{2\alpha}N< \frac{2\alpha}{N-2\alpha}$. Therefore Theorem \ref{teo 1.1} does not provide any description of $u_\infty$ in the region \begin{center}

${\cal U}:=\left\{(\alpha,p)\in (0,1)\times (1,\frac{N}{N-2}):\frac N{2\alpha}<1+\frac{2\alpha}N,\frac {N}{2\alpha}<p<1+\frac{2\alpha}{N}\right\}.
$
\end{center}
Furthermore, in  parts $(ii)$ and $(iii)$, we do not obtain that $u_\infty=u_s$, since  (\ref{1.2.1}) and (\ref{1.1.1}) do not provide sharp estimates on $u_\infty$ in order it to belong to the uniqueness class  characterized by \re{strong singular unique}. \medskip

The paper is organized as follows. In Section 2, we present some  some estimates
for the Green kernel and comparison principles. In Section 3, we prove that the  weak solution of (\ref{eq1.1})
is a classical solution of (\ref{eq1.2}).
Section  4 is devoted to analyze the limit of weakly singular  solutions as $k\to\infty$.

\setcounter{equation}{0}
\section{Preliminaries}
The purpose of this section is to recall some known results. We denote by $B_r(x)$ the
ball centered at $x$ with radius $r$ and $B_r:=B_r(0)$.

\begin{lemma}\label{lemma 3} Assume that  $0<p<p^*_\alpha,$
then there exists $c_4,c_5,c_6>1$  such that\\
$(i)$ if $p\in(0,\frac{2\alpha}{N-2\alpha})$,
$$
\frac1{c_4}\le \mathbb{G}_\alpha[(\mathbb{G}_\alpha[\delta_0])^p]\le c_4\quad \ {\rm in}\ \
B_r\setminus\{0\};
$$
$(ii)$ if $p=\frac{2\alpha}{N-2\alpha}$,
$$
-\frac1{c_5}\ln|x|\le \mathbb{G}_\alpha[(\mathbb{G}_\alpha[\delta_0])^p]\le -{c_5}\ln|x|\quad \ {\rm in}\ \ B_r\setminus\{0\};
$$
$(iii)$ if $p\in(\frac{2\alpha}{N-2\alpha}, p^*_\alpha)$,
$$
\frac1{c_6}|x|^{2\alpha-(N-2\alpha)p}\le
\mathbb{G}_\alpha[(\mathbb{G}_\alpha[\delta_0])^p]\le
{c_6}|x|^{2\alpha-(N-2\alpha)p}\quad \ {\rm in}\ \ B_r\setminus\{0\},
$$
where $r=\frac14\min\{1,dist(0,\partial\Omega)\}$ and $\mathbb{G}_\alpha$ is defined by (\ref{1.3}).
\end{lemma}
{\bf Proof.} The proof follows easily from Chen-Song's estimates of Green functions \cite{CS}, see  \cite[Theorem 5.2]{CV0} for a detailled computation.\hfill$\Box$

\begin{teo}\label{teo CP}
Assume that $O$ is a bounded domain of $\R^N$ and
 $u_1$, $u_2$ are continuous in $\bar O$ and satisfy
$$(-\Delta)^\alpha  u+u^p=0 \quad{\rm in}\quad O.$$
Moreover, we assume that $u_1\ge u_2$ in $O^c$.
Then, \smallskip

\noindent (i) either $u_1> u_2\quad{\rm in}\quad O,$\smallskip

\noindent (ii) or $u_1\equiv u_2\quad{\rm\ a.e.\ in}\quad \R^N.$
\end{teo}
{\bf Proof.} The proof refers to \cite[Theorem 2.3]{CFQ} (see also  \cite[Theorem 5.2]{CS1}).
\hfill$\Box$\smallskip

The following stability result is proved in \cite[Theorem 2.2]{CFQ}.

\begin{teo}\label{stability}
Suppose that $\mathcal{O}$ is a bounded $C^2$ domain and
$h:\R\to \R$ is continuous.  Assume $\{u_n\}$ is a sequence of functions, uniformly  bounded in
$L^1(\mathcal{O}^c,\frac{dy}{1+|y|^{N+2\alpha}})$, satisfying
$$(-\Delta)^\alpha u_n+ h(u_n)\ge f_n\; (\mbox{resp } (-\Delta)^\alpha
u_n+ h(u_n)\le f_n\;)\quad\mbox{in }\mathcal{O}$$
in the viscosity sense, where the $f_n$
are continuous in $\mathcal{O}$. If there holds

\noindent (i) $u_n\to u$ locally uniformly in  $\mathcal{O}$,

\noindent  (ii) $ u_n\to u$  in  $L^1(\R^N,\frac{dy}{1+|y|^{N+2\alpha}})$,

\noindent (iii) $ f_n\to f$ locally uniformly in  $\mathcal{O}$,\smallskip

\noindent then
    $$(-\Delta)^\alpha u+ h(u)\ge f\; (\mbox{resp } (-\Delta)^\alpha
u+ h(u)\le f\;)\quad\mbox{in }\mathcal{O}$$
in the viscosity sense.
\end{teo}

\setcounter{equation}{0}
\section{Regularity}

In this section, we  prove that any weak solution of (\ref{eq1.1}) is a classical solution of (\ref{eq1.2}).
To this end, we introduce some auxiliary lemma.
\begin{lemma}\label{lm 2.1}
Assume that $w\in C^{2\alpha+\epsilon}(\bar B_1)$ with $\epsilon>0$
 satisfies
 $$(-\Delta)^\alpha w=h\quad {\rm in}\quad B_1,$$
 where $h\in C^1(\bar B_1)$. Then for $\beta\in (0,2\alpha)$, there exists ${c_7}>0$ such that
\begin{equation}\label{2.0}
\|w\|_{C^\beta(\bar B_{1/4})}\le {c_7}(\|w\|_{L^\infty(B_1)}+\|h\|_{L^\infty(B_1)}+\|(1+|\cdot|)^{-N-2\alpha}w\|_{L^1(\R^N)}).
\end{equation}

\end{lemma}
{\bf Proof.}
Let $\eta:\R^N\to[0,1]$ be a $C^\infty$  function such that
$$\eta=1\quad {\rm in}\quad  B_{\frac34}\quad {\rm and}\quad   \eta=0\quad {\rm in}\quad  B_1^c. $$
We denote $v=w\eta$, then $v\in C^{2\alpha+\epsilon}(\R^N)$ and for $x\in B_{\frac12}$, $\epsilon\in(0,\frac14)$,
\begin{eqnarray*}
  (-\Delta)_\epsilon^\alpha v(x) &=& -\int_{\R^N\setminus B_\epsilon}\frac{v(x+y)-v(x)}{|y|^{N+2\alpha}}dy
    \\   &=&(-\Delta)_\epsilon^\alpha w(x)+\int_{\R^N\setminus B_\epsilon}\frac{(1-\eta(x+y))w(x+y)}{|y|^{N+2\alpha}}dy.
 \end{eqnarray*}
Together with the fact of $\eta(x+y)=1$ for $y\in B_\epsilon$, we have
 $$\int_{\R^N\setminus B_\epsilon}\frac{(1-\eta(x+y))w(x+y)}{|y|^{N+2\alpha}}dy=\int_{\R^N}\frac{(1-\eta(x+y))w(x+y)}{|y|^{N+2\alpha}}dy=: h_1(x),$$
thus, $$(-\Delta)^\alpha v=h+h_1\quad{\rm in}\quad B_{\frac12}.$$
For $x\in B_{\frac12}$ and $z\in\R^N\setminus B_{\frac34}$, there holds
\begin{eqnarray*}
  |z-x|\ge |z|-|x| \ge |z|-\frac12
   \ge \frac1{16}(1+|z|)
\end{eqnarray*}
which implies
\begin{eqnarray*}
  |h_1(x)| = \abs{\int_{\R^N}\frac{(1-\eta(z))w(z)}{|z-x|^{N+2\alpha}}dz}
   &\le& \int_{\R^N\setminus B_{\frac34}}\frac{|w(z)|}{|z-x|^{N+2\alpha}}dz \\
   &\le & 16^{N+2\alpha} \int_{\R^N}\frac{|w(z)|}{(1+|z|)^{N+2\alpha}}dz\\
   &=&16^{N+2\alpha}\|(1+|\cdot|)^{-N-2\alpha}w\|_{L^1(\R^N)}.
\end{eqnarray*}
By \cite[Proposition 2.1.9]{S},  for $\beta\in (0,2\alpha)$, there exists ${c_8}>0$ such that
\begin{eqnarray*}
\|v\|_{C^\beta(\bar B_{1/4})}&\le& {c_8}(\|v\|_{L^\infty(\R^N)}+\|h+h_1\|_{L^\infty(B_{1/2})}) \\
   &\le& {c_8}(\|w\|_{L^\infty(B_1)}+\|h\|_{L^\infty(B_1)}+\|h_1\|_{L^\infty(B_{1/2})})\\
   &\le& {c_9}(\|w\|_{L^\infty(B_1)}+\|h\|_{L^\infty(B_1)}+\|(1+|\cdot|)^{-N-2\alpha}w\|_{L^1(\R^N)}),
\end{eqnarray*}
where  $ {c_9}=16^{N+2\alpha}{c_8}$.
Combining with $w=v$ in $B_{\frac34}$,  we obtain (\ref{2.0}). \hfill$\Box$

\begin{teo}\label{pr 2.1}
Let $\alpha\in(0,1)$ and  $0<p<p^*_\alpha$,
then the weak solution of (\ref{eq1.1}) is a classical solution of (\ref{eq1.2}).
\end{teo}
{\bf Proof.} Let $u_k$ be the weak solution of (\ref{eq1.1}). By \cite[Theorem 1.1]{CV2},
we have
\begin{equation}\label{2.1}
0\le u_k= \mathbb{G}_\alpha[k\delta_0]-\mathbb{G}_\alpha[u_k^p ]\le\mathbb{G}_\alpha[k\delta_0].
\end{equation}
We observe that
$\mathbb{G}_\alpha[k\delta_0]=k\mathbb{G}_\alpha[\delta_0]=kG_\alpha(\cdot,0)$ is $C^2_{loc}(\Omega\setminus\{0\})$.
Denote by $O$ an open set satisfying $ \bar O\subset \Omega\setminus B_r$ with $r>0$.
Then $\mathbb{G}_\alpha[k\delta_0]$  is uniformly bounded in $\Omega\setminus B_{r/2}$,  so is $u_k^p$ by (\ref{2.1}).

 Let $\{g_n\}$ be a sequence nonnegative functions in $C^\infty_0(\R^N)$  such that
$g_n\to\delta_0$ in the weak sense of measures and  let $w_n$ be the solution of
\begin{equation}
 \arraycolsep=1pt
\begin{array}{lll}
 (-\Delta)^\alpha  u+u^p=kg_n\quad & \rm{in}\quad\Omega\\[2mm]
 \phantom{   (-\Delta)^\alpha  + u^p}
u=0\quad & \rm{in}\quad \Omega^c.
\end{array}
\end{equation}
From \cite{CV2}, we obtain that
\begin{equation}\label{2.2}
u_k=\lim_{n\to\infty} w_n\quad {\rm a.e.\ in}\quad \Omega.
\end{equation}
 We observe that $0\le w_n= \mathbb{G}_\alpha[kg_n]-\mathbb{G}_\alpha[w_n^p ]\le k\mathbb{G}_\alpha[g_n]$ and
 $\mathbb{G}_\alpha[g_n]$ converges to $\mathbb{G}_\alpha[\delta_0]$ uniformly in any compact set of $\Omega\setminus \{0\}$ and in $L^1(\Omega)$; then there exists $c_{10}>0$ independent of $n$ such that
$$\norm{w_n}_{L^\infty(\Omega\setminus B_{r/2})}\le {c_{10}}k\quad{\rm and}\quad \norm{w_n}_{L^1(\Omega)}\le {c_{10}}k.$$
By \cite[Corollary 2.4]{RS} and Lemma \ref{lm 2.1}, there exist $\epsilon>0$, $\beta\in(0,2\alpha)$ and positive constants ${c_{11}},{c_{12}},{c_{13}}>0$ independent of $n$ and $k$, such that
\begin{eqnarray*}
 &&\norm{w_n}_{C^{2\alpha+\epsilon}(O)} \le {c_{11}}(\norm{w_n}^p_{L^\infty(\Omega\setminus B_{\frac r2})}+\norm{kg_n}_{L^{\infty}(\Omega\setminus B_{\frac r2})}+\|w_n\|_{C^\beta(\Omega\setminus B_{\frac {3r}4})})\\
 && \quad\le {c_{12}}(\norm{w_n}^p_{L^\infty(\Omega\setminus B_{\frac r2})}+\norm{w_n}_{L^\infty(\Omega\setminus B_{\frac r2})}+\norm{kg_n}_{L^{\infty}(\Omega\setminus B_{\frac r2})}+\|w_n\|_{L^1(\Omega)})\\
  &&\quad\le {c_{13}}(k+k^p).
\end{eqnarray*}
Therefore,  together with (\ref{2.2}) and the Arzela-Ascoli Theorem, it follows that
$u_k\in C^{2\alpha+\frac\epsilon2}(O)$. This implies that $u_k$ is $C^{2\alpha+\frac\epsilon2}$ locally in $\Omega\setminus\{0\}$.
Therefore, $w_n\to u_k$ and $g_n\to 0$ uniformly in any compact subset of $\Omega\setminus\{0\}$ as $n\to\infty$.
We conclude  that $u_k$ is a classical solution of (\ref{eq1.2}) by Theorem \ref{stability}.
\hfill$\Box$
\begin{corollary}\label{cr 1}
Let $u_k$ be the weak solution of (\ref{eq1.1}) and  $O$ be an open set satisfying $ \bar O\subset \Omega\setminus B_r$ with $r>0$. Then there exist $\epsilon>0$ and ${c_{14}}>0$ independent of $k$ such that
  \begin{equation}\label{hu 2}
 \norm{u_k}_{C^{2\alpha+\epsilon}(O)} \le {c_{14}}(\norm{u_k}^p_{L^\infty(\Omega\setminus B_{\frac r2})}+\norm{u_k}_{L^\infty(\Omega\setminus B_{\frac r2})}+\|u_k\|_{L^1(\Omega)}).
\end{equation}
\end{corollary}
{\bf Proof.}  By Theorem \ref{pr 2.1}, $u_k$ is a solution of (\ref{eq1.2}).
Then the result follows from \cite[Corollary 2.4]{RS} and Lemma \ref{lm 2.1} since there exist $\epsilon>0$, $\beta\in(0,2\alpha)$ and constants ${c_{15}},{c_{16}}>0$, independent of $k$, such that
\begin{eqnarray*}
\norm{u_k}_{C^{2\alpha+\epsilon}(O)} &\le& {c_{15}}(\norm{u_k}^p_{L^\infty(\Omega\setminus B_{\frac r2})}+\|u_k\|_{C^\beta(\Omega\setminus B_{\frac {3r}4})})\\
 &\le& {c_{16}}(\norm{u_k}^p_{L^\infty(\Omega\setminus B_{\frac r2})}+\norm{u_k}_{L^\infty(\Omega\setminus B_{\frac r2})}+\|u_k\|_{L^1(\Omega)}).
\end{eqnarray*}
  \hfill$\Box$

\begin{teo}\label{teo 2}
Assume that the weak solutions $u_k$ of (\ref{eq1.1}) satisfy
\begin{equation}\label{ying 1}
\norm{u_k}_{L^1(\Omega)}\le {c_{17}}
\end{equation}
for some ${c_{17}}>0$ independent of $k$ and that for any $r\in(0,dist(0,\partial\Omega))$, there exists ${c_{18}}>0$ independent of $k$ such that
\begin{equation}\label{ying 2}
\norm{u_k}_{L^\infty(\Omega\setminus B_{\frac{r}{2}})}\le {c_{18}}.
\end{equation}
Then $u_\infty$ is a classical solution of (\ref{eq1.2}).
\end{teo}
{\bf Proof.} Let $O$ be an open set satisfying $ \bar O\subset \Omega\setminus B_{r}$ for $0<r<dist(0,\partial\Omega)$.
By (\ref{hu 2}), (\ref{ying 1}) and (\ref{ying 2}),
there  exist $\epsilon>0$ and ${c_{19}}>0$ independent of $k$ such that
$$\norm{u_k}_{C^{2\alpha+\epsilon}(O)}\le {c_{19}}.$$
Together with (\ref{definition infty}) and the Arzela-Ascoli Theorem, it implies that
$u_\infty$ belongs to $C^{2\alpha+\frac\epsilon2}(O)$. Hence $u_\infty$ is $C^{2\alpha+\frac\epsilon2}$, locally in $\Omega\setminus\{0\}$.
Therefore, $w_n\to u_k$ and $g_n\to 0$ uniformly in any compact set of $\Omega\setminus\{0\}$ as $n\to\infty$.
Applying Theorem \ref{stability} we conclude that $u_\infty$ is a classical solution of (\ref{eq1.2}).
\hfill$\Box$
\setcounter{equation}{0}
\section{The limit of weakly singular solutions}
We recall that $u_k$  denotes the weak solution of (\ref{eq1.1}) and $d=\min\{1,dist(0,\partial\Omega)\}.$

\subsection{The case $p\in(0,1+\frac{2\alpha}N]$}

\begin{proposition}\label{pr 3.1}
Let $p\in(0,1]$, then $\lim_{k\to\infty}u_k(x)=\infty$ for $x\in \Omega.$
\end{proposition}
{\bf Proof.}  We observe that $\mathbb{G}_\alpha[\delta_0], \mathbb{G}_\alpha[(\mathbb{G}_\alpha[\delta_0])^p]>0$ in $\Omega$.
Since by (\ref{1.3}) 
\begin{eqnarray*}
u_k\ge k\mathbb{G}_\alpha[\delta_0]- k^p\mathbb{G}_\alpha[(\mathbb{G}_\alpha[\delta_0])^p],
\end{eqnarray*}
this implies the claim when  $p\in(0,1)$, for any $x\in  \Omega$.
For $p=1$, $u_k=ku_1$. The proof follows since $u_1>0$ in $\Omega$. \hfill$\Box$

\bigskip
Now we consider  the case of $p\in(1,1+\frac{2\alpha}N]$.  Let $\{r_k\}\subset(0,\frac d2]$ be a strictly decreasing sequence of  numbers satisfying $\lim_{k\to\infty}r_k=0$.
Denote by $\{z_k\}$  the sequence of functions defined by
\begin{equation}\label{4.2}
\arraycolsep=1pt
z_{k}(x)=\left\{
\begin{array}{lll}
-d^{-N},\quad & x\in B_{r_k}\\[2mm]
|x|^{-N}-d^{-N},\quad & x\in B_{r_k}^c.
\end{array}
\right.
\end{equation}

\begin{lemma}\label{lm 3.1}
Let $\{\rho_k\}$ be a strictly decreasing sequence of  numbers  such that $\frac{r_k}{\rho_k}<\frac12$ and $\lim_{k\to\infty}\frac{r_k}{\rho_k}=0$.
Then
$$(-\Delta)^\alpha z_{k}(x)\le -c_{1,k}|x|^{-N-2\alpha},\quad\quad x\in B_{\rho_k}^c$$
where $c_{1,k}= -c_{20}\log(\frac{r_k}{\rho_k})$ with $c_{20}>0$ independent of $k$.
\end{lemma}
{\bf Proof.} For any $x\in B_{\rho_k}^c$, there holds
\begin{eqnarray*}
 &&(-\Delta)^\alpha z_{k}(x)= -\frac12\int_{\R^N}\frac{z_k(x+y)+z_k(x-y)-2z_k(x)}{|y|^{N+2\alpha}}dy  \\
 &&\qquad= -\frac12\int_{\R^N}\frac{|x+y|^{-N}\chi_{B^c_{r_k}(-x)}(y)+|x-y|^{-N}\chi_{B^c_{r_k}(x)}(y)-2|x|^{-N}}{|y|^{N+2\alpha}}dy \\
\\&&\qquad=-\frac12|x|^{-N-2\alpha}\int_{\R^N}\frac{\delta(x,z,r_k)}{|z|^{N+2\alpha}}dz,
\end{eqnarray*}
where $\delta(x,z,r_k)=|z+e_x|^{-N}\chi_{B^c_{\frac{r_k}{|x|}}(-e_x)}(z)+|z-e_x|^{-N}\chi_{B^c_{\frac{r_k}{|x|}}(e_x)}(z)-2$ and $e_x=\frac x{|x|}$.

We observe that $\frac{r_k}{|x|}\le\frac{r_k}{\rho_k} <\frac12 $ and $|z\pm e_x|\ge 1-|z|\ge \frac12$ for $z\in B_{\frac12}$. Then there exists $c_{21}>0$ such that
\begin{eqnarray*}
|\delta(x,z,r_k)|=||z+e_x|^{-N}+|z-e_x|^{-N}-2|
\le c_{21}|z|^2.
\end{eqnarray*}
Therefore,
\begin{eqnarray*}
|\int_{B_{\frac12}(0)}\frac{\delta(x,z,r_k)}{|z|^{N+2\alpha}}dz | &\le &\int_{B_{\frac12}(0)}\frac{|\delta(x,z,r_k)|}{|z|^{N+2\alpha}}dz
\\&\le& c_{21}\int_{B_{\frac12}(0)} |z|^{2-N-2\alpha} dz\le c_{22},
\end{eqnarray*}
where  $c_{22}>0$ is independent of  $k$.

\noindent When $z\in B_{\frac12}(-e_x)$ there holds
\begin{eqnarray*}
\int_{B_{\frac12}(-e_x)}\frac{\delta(x,z,r_k)}{|z|^{N+2\alpha}}dz  &\ge &\int_{B^c_{\frac12}(-e_x)}\frac{|z+e_x|^{-N}\chi_{B_{\frac{r_k}{|x|}}(-e_x)}(z)-2} {|z|^{N+2\alpha}}dz
\\&\ge&c_{23} \int_{B_{\frac12}(0)\setminus B_{\frac{r_k}{|x|}}(0)}(|z|^{-N}-2)dz
\\&\ge&-c_{24} \log(\frac{r_k}{|x|})\ge -c_{24} \log(\frac{r_k}{\rho_k}),
\end{eqnarray*}
where  $c_{23},c_{24}>0$ are independent of  $k$.

\noindent For $z\in B_{\frac12}(e_x)$, we have
$$\int_{B_{\frac12}(e_x)}\frac{\delta(x,z,r_k)}{|z|^{N+2\alpha}}dz=\int_{B_{\frac12}(-e_x)}\frac{\delta(x,z,r_k)}{|z|^{N+2\alpha}}dz.$$

\noindent Finally, for $z\in O:=\R^N\setminus (B_{\frac12}(0)\cup B_{\frac12}(-e_x)\cup B_{\frac12}(e_x))$,
we have
\begin{eqnarray*}
|\int_{O}\frac{\delta(x,z,r_k)}{|z|^{N+2\alpha}}dz| \le  c_{25} \int_{B_{\frac12}^c(0)}\frac{|z|^{-N}+1}{|z|^{N+2\alpha}}dz
\le c_{26},
\end{eqnarray*}
where $c_{25},c_{26}>0$ are independent of $k$.

\noindent Combining these inequalities we obtain that there exists $c_{20}>0$ independent of $k$ such that
$$(-\Delta)^\alpha z_{k}(x)|x|^{N+2\alpha}\le c_{20}\log(\frac{r_k}{\rho_k}):=c_{1,k},$$
which ends the proof.\hfill$\Box$

\begin{proposition}\label{pr 4.0}
Assume that
\begin{equation}\label{hu 1.1}
\frac{2\alpha}{N-2\alpha}<1+\frac{2\alpha}N,\quad  \max\{1,\frac{2\alpha}{N-2\alpha}\}<p<1+\frac{2\alpha}N
\end{equation}
 and $z_k$ is defined by (\ref{4.2}) with $r_k=k^{-\frac{p-1}{N-(N-2\alpha)p}}(\log k)^{-2} $.
Then there exists $k_0>0$ such that for any $k\ge k_0$
\begin{equation}\label{E1}u_k\ge c_{2,k}^{\frac1{p-1}}z_k\quad \ {\rm in}\quad B_d,\end{equation}
where $c_{2,k}=\ln\ln k$.
\end{proposition}
{\bf Proof.}
For  $p\in(\max\{1,\frac{2\alpha}{N-2\alpha}\},1+\frac{2\alpha}N)$, it follows by
 (\ref{1.3})  and Lemma \ref{lemma 3}-$(iii)$ that there exist $\rho_0\in(0,d)$ and $c_{27},c_{28}>0$ independent of $k$ such that,  for  $x\in \bar B_{\rho_0}\setminus\{0\}$,
\begin{eqnarray*}
 u_k(x)&\ge& k\mathbb{G}_\alpha[\delta_0](x)-k^p\mathbb{G}_\alpha[(\mathbb{G}_\alpha[\delta_0])^p](x) \\
 &\ge& c_{27}k|x|^{-N+2\alpha}-c_{28}k^p|x|^{-(N-2\alpha)p+2\alpha}  \\
   &=& c_{27}k|x|^{-N+2\alpha}(1-\frac{c_{28}}{c_{27}}k^{p-1}|x|^{N-(N-2\alpha)p}).
\end{eqnarray*}
We choose
\begin{equation}\label{hu 1}
\rho_k=k^{-\frac{p-1}{N-(N-2\alpha)p}}(\log k)^{-1}.
\end{equation}
There exits $k_1>1$ such that for $k\ge k_1$
\begin{eqnarray}
 u_k(x)&\ge& c_{27}k|x|^{-N+2\alpha}(1-\frac{c_{28}}{c_{27}}k^{p-1}\rho_k^{N-(N-2\alpha)p})\nonumber
  \\&\ge&\frac{c_{27}}2k|x|^{-N+2\alpha},\quad x\in \bar B_{\rho_k}\setminus\{0\}. \label{chen 1}
\end{eqnarray}
Since $p<1+\frac{2\alpha}{N}$, $1-\frac{2\alpha(p-1)}{N-(N-2\alpha)p}>0$ and there exists $k_0\ge k_1$ such that
\begin{equation}\label{chen 2}
\frac{c_{27}}2kr_k^{2\alpha}\ge (\ln\ln k)^{\frac1{p-1}},
\end{equation}
for $k\geq k_0$. This implies
$$\frac{c_{27}}2k|x|^{2\alpha}\ge (\ln\ln k)^{\frac1{p-1}},\quad x\in \bar B_{\rho_k}\setminus B_{r_k}.$$
Together with  (\ref{4.2}) and (\ref{chen 1}), we derive
$$u_k(x)\ge (\ln\ln k)^{\frac1{p-1}}z_k(x),\quad x\in \bar B_{\rho_k}\setminus B_{r_k},$$
for $k\ge k_0$.
Furthermore, it is clear that 
$$(\ln\ln k)^{\frac1{p-1}}z_k(x)\le0\le u_k(x)$$
whenever $x\in B_{r_k}$ or $x\in B_d^c$. Set $c_{2,k}=\ln\ln k$, then by Lemma \ref{lm 3.1}
\begin{eqnarray*}
  (-\Delta)^\alpha c_{2,k}^{\frac1{p-1}}z_k(x)+c_{2,k}^{\frac p{p-1}}z_k(x)^p \le c_{2,k}^{\frac p{p-1}} |x|^{-N-2\alpha}
  (-1+|x|^{N+2\alpha-Np})\le 0,
\end{eqnarray*}
for any $x\in B_d\setminus B_{\rho_k}$, since $N+2\alpha-Np\ge 0$ and $d\leq 1$.
Applying Theorem \ref{teo CP}, we infer that
$$c_{2,k}^{\frac1{p-1}}z_k(x)\le u_k(x)\qquad \forall x\in\bar B_d,$$
which ends the proof.
\hfill$\Box$

\begin{proposition}\label{pr 4.1}
Assume 
\begin{equation}\label{hu 1.2}
  1<\frac{2\alpha}{N-2\alpha}\le1+\frac{2\alpha}N\,\mbox{ and }\; p=\frac{2\alpha}{N-2\alpha}
\end{equation}
and let $z_k$ be defined by (\ref{4.2}) with $r_k= k^{-\frac{2\alpha}{N(N-2\alpha)}}(\log k)^{-3}$ and $k>2$.
Then there exists $k_0>2$ such that
(\ref{E1}) holds for $k\ge k_0$.
\end{proposition}
{\bf Proof.}
 By (\ref{1.3}) and Lemma \ref{lemma 3}-$(ii)$,  there exist $\rho_0\in(0,d)$ and $c_{30},c_{31}>0$ independent of $k$, such that for  $x\in \bar B_{\rho_0}\setminus\{0\}$
\begin{eqnarray*}
 u_k(x)&\ge& c_{30}k|x|^{-N+2\alpha}+c_{31}k^p\log|x| \\
   &=& c_{30}k|x|^{-N+2\alpha}(1+\frac{c_{31}}{c_{30}}k^{p-1}|x|^{N-2\alpha}\log|x|).
\end{eqnarray*}
If we choose  $\rho_k=k^{-\frac{2\alpha}{N(N-2\alpha)}}(\log k)^{-2}$ there exists $k_1>1$ such that for $k\ge k_1$, we have
$1+\frac{c_{31}}{c_{30}}k^{p-1}\rho_k^{N-2\alpha}\log(\rho_k)\ge \frac12$
and
\begin{equation}\label{chen 01}
 u_k(x) \ge\frac{c_{30}}2k|x|^{-N+2\alpha}\qquad \forall x\in \bar B_{\rho_k}\setminus\{0\}.
\end{equation}
Since $\frac{2\alpha}{N-2\alpha}<1+\frac{2\alpha}N $, there holds $1-\frac{4\alpha^2}{N(N-2\alpha)}>0$ and there exists $k_0\geq k_1$ such that
$$\frac{c_{30}}2kr_k^{2\alpha}=\frac{c_{30}}2k^{1-\frac{4\alpha^2}{N(N-2\alpha)}}(\log k)^{-6\alpha}\geq (\ln\ln k)^{\frac{1}{p-1}}$$
 for $k\ge k_0$.
The remaining of the proof is the same as in Proposition \ref{pr 4.0}.
\hfill$\Box$
\medskip


In the sequel, we point out the fact that  the limit behavior of the $u_k$ depends which of the following three cases holds:
\begin{equation}\label{chen 1.1}
\frac{2\alpha}{N-2\alpha}=1+\frac{2\alpha}N=\frac N{2\alpha};
\end{equation}
\begin{equation}\label{chen 1.2}
\frac{2\alpha}{N-2\alpha}<1+\frac{2\alpha}N<\frac N{2\alpha};
\end{equation}
\begin{equation}\label{yan 0}
 \frac{2\alpha}{N-2\alpha}>1+\frac{2\alpha}N>\frac N{2\alpha}.
\end{equation}

\begin{proposition}\label{pr 4.2}
Assume
\begin{equation}\label{hu 1.3}
  1<\frac{2\alpha}{N-2\alpha}\le 1+\frac{2\alpha}N\,\mbox{ and }\; 1<p< \frac{2\alpha}{N-2\alpha},
\end{equation}
 or
 \begin{equation}\label{hu 1.4}
1+\frac{2\alpha}N< \frac{2\alpha}{N-2\alpha}\,\mbox{ and }\; 1<p< \frac N{2\alpha},
\end{equation}
and $z_k$ is defined by (\ref{4.2}) with $r_k= k^{-\frac{p-1}{N-2\alpha}}(\log k)^{-1}$.
 Then there exists $k_0>2$ such that
(\ref{E1}) holds for $k\ge k_0$.
\end{proposition}
{\bf Proof.}
 By (\ref{1.3}) and Lemma \ref{lemma 3}-$(i)$,  there exist $\rho_0\in(0,d)$ and $c_{33},c_{34}>0$ independent of $k$ such that for  $x\in \bar B_{\rho_0}\setminus\{0\}$,
\begin{eqnarray*}
 u_k(x)&\ge& c_{33}k|x|^{-N+2\alpha}-c_{34}k^p \\
   &=& c_{33}k|x|^{-N+2\alpha}(1-\frac{c_{34}}{c_{33}}k^{p-1}|x|^{N-2\alpha}).
\end{eqnarray*}
We choose  $\rho_k=k^{-\frac{p-1}{N-2\alpha}}$. Then there exists $k_1>1$ such that for $k\ge k_1$,
$1-\frac{c_{34}}{c_{33}}k^{p-1}\rho_k^{N-2\alpha}\ge \frac12$
and
\begin{equation}\label{chen 001}
 u_k(x) \ge\frac{c_{33}}2k|x|^{-N+2\alpha}\qquad \forall x\in \bar B_{\rho_k}\setminus\{0\}.
\end{equation}
Clearly $p<\frac N{2\alpha} $ by assumptions (\ref{hu 1.3}), (\ref{hu 1.4}), together with relations (\ref{chen 1.1})(\ref{chen 1.2}), (\ref{yan 0}), thus $1-(p-1)\frac{2\alpha}{N-2\alpha}>0$. Therefore  there exists $k_0\ge k_1$ such that
$$\frac{c_{33}}2kr_k^{2\alpha}=\frac{c_{33}}2k^{1-(p-1)\frac{2\alpha}{N-2\alpha}}(\log k)^{-2\alpha}\ge (\log\log k)^{\frac1{p-1}}=c_{2,k}^{\frac{1}{p-1}}$$
for $k\geq k_0$. The remaining of the  proof is similar to the one of Proposition \ref{pr 4.0}.
\hfill$\Box$
\medskip

\subsection{The case $p\in(1,1+\frac{2\alpha}N]$}

We give below the proof, in two steps, of Theorem \ref{teo 1}  part $(i)$ with $p\in(1,1+\frac{2\alpha}N]$ and Theorem \ref{teo 1.1} part $(i)$ with $p\in(1,\frac{2\alpha}N]$.
\smallskip

\noindent{\it Step 1: We claim that $u_{\infty}=\infty$ in $B_d$}. 
We observe that for $\frac{2\alpha}{N-2\alpha}< 1+\frac{2\alpha}N$, Propositions \ref{pr 4.0}, \ref{pr 4.1}, \ref{pr 4.2} cover the case $p\in(\max\{1,\frac{2\alpha}{N-2\alpha}\},1+\frac{2\alpha}N)$, the case $1<\frac{2\alpha}{N-2\alpha}<1+\frac{2\alpha}N$ along with  $p=\frac{2\alpha}{N-2\alpha}$ and the case $1<\frac{2\alpha}{N-2\alpha}<1+\frac{2\alpha}N$ along with $p\in(1, \frac{2\alpha}{N-2\alpha})$ respectively. For  $\frac{2\alpha}{N-2\alpha}= 1+\frac{2\alpha}N$,  Proposition \ref{pr 4.1}, \ref{pr 4.2} cover the case  $p=\frac{2\alpha}{N-2\alpha}$
and the case $p\in(1, \frac{2\alpha}{N-2\alpha})$ respectively. So it covers $p\in(1,1+\frac{2\alpha}N]$ in Theorem \ref{teo 1}  part $(i)$.  When $\frac{2\alpha}{N-2\alpha}> 1+\frac{2\alpha}N$, Proposition \ref{pr 4.2} covers $p\in(1, \frac N{2\alpha})$ in  Theorem \ref{teo 1.1} part $(i)$. Therefore,  we have
$$u_{\infty}\ge  c_{2,k}^{\frac1{p-1}}z_k\quad {\rm in}\quad B_d$$
and since for any $x\in B_d\setminus\{0\}$, $\lim_{k\to\infty}c_{2,k}^{\frac1{p-1}}z_k(x)=\infty,$
we deive
$$u_{\infty}=\infty \quad {\rm in}\quad B_d.$$

\noindent {\it Step 2: We claim that $u_{\infty}=\infty$ in $\Omega$}. By the  fact of $u_{\infty}=\infty$ in $B_d$ and $u_{k+1}\ge u_k$ in $\Omega$,
then for any $n>1$
there exists $k_n>0$ such that $u_{k_n}\ge n$ in $B_d$.
For any $x_0\in \Omega\setminus  B_d$, there exists $\rho>0$ such that
$\bar B_\rho(x_0)\subset \Omega\cap  B_{d/2}^c$.
 We denote by $w_n$ the solution of
\begin{equation}\label{5.1}
\arraycolsep=1pt
\begin{array}{lll}
 (-\Delta)^\alpha  u+u^p=0 \quad & {\rm in}\quad  B_{\rho}(x_0)\\[2mm]
 \phantom{  (-\Delta)^\alpha  +u^p}
u=0  \quad & {\rm in}\quad B_{\rho}^c(x_0)\setminus B_{d/2}\\[2mm]
\phantom{,  (-\Delta)^\alpha  +u^p}
u=n  \quad & {\rm in}\quad B_{d/2}.
\end{array}
\end{equation}
Then by Theorem \ref{teo CP}, we have
\begin{equation}\label{4.1.3}
 u_{k_n}\ge w_n.
\end{equation}
Let $\eta_1$ be the solution of
$$
\arraycolsep=1pt
\begin{array}{lll}
 (-\Delta)^\alpha  u=1 \quad & {\rm in}\quad  B_{\rho}(x_0)\\[2mm]
 \phantom{  (-\Delta)^\alpha  }
u=0  \quad & {\rm in}\quad  B^c_{\rho}(x_0),
\end{array}
$$
and $v_n=w_n-n\chi_{B_{d/2}},$
then $v_n=w_n$ in $B_\rho(x_0)$ and
\begin{eqnarray*}
      (-\Delta)^\alpha v_n(x)+v_n^p(x) &=&  (-\Delta)^\alpha w_n(x)- n(-\Delta)^\alpha \chi_{B_{d/2}}(x) +w_n^p(x) \\
        &=&n\int_{B_{d/2}}\frac{dy}{|y-x|^{N+2\alpha}}\qquad \forall x\in B_\rho(x_0).
     \end{eqnarray*}
This means that  $v_n$  is a solution of
\begin{equation}\label{4.1.2}
\arraycolsep=1pt
\begin{array}{lll}
\displaystyle (-\Delta)^\alpha  u+u^p=n\int_{B_{d/2}}\frac{dy}{|y-x|^{N+2\alpha}} \quad & {\rm in}\quad  B_{\rho}(x_0),\\[2mm]
 \phantom{  (-\Delta)^\alpha  +u^{p,}}
u=0  \quad & {\rm in}\quad B_{\rho}^c(x_0).
\end{array}
\end{equation}
It is clear that 
$$\frac1{c_{35}}\le \int_{B_{d/2}}\frac{dy}{|y-x|^{N+2\alpha}}\le c_{35}\qquad\forall x\in B_{\rho}(x_0)$$ 
for some $c_{35}>1$. Furhermore
 $(\frac{n}{2c_{35}\max\eta_1})^{\frac1p} \eta_1$ is sub solution of (\ref{4.1.2}) for $n $ large enough. Then using Theorem \ref{teo CP}, we obtain that
 $$v_n\ge (\frac{n}{2c_{35}\max\eta_1})^{\frac1p} \eta_1\qquad \forall  x\in B_{\rho}(x_0),$$
 which implies that
 $$w_n\ge (\frac{n}{2c_{35}\max\eta_1})^{\frac1p} \eta_1\qquad \forall  x\in B_{\rho}(x_0).$$ 
 Then
 $$\lim_{n\to\infty}w_n(x_0)\to\infty.$$
Since $x_0$ is arbitrary and together with (\ref{4.1.3}), it implies that $u_\infty=\infty$ in $\Omega$, which completes the proof.
\hfill$\Box$

\subsection{The case of $p\in(1+\frac{2\alpha}N, \frac{N}{N-2\alpha})$}

\begin{proposition}\label{pr 3.2}
Let $\alpha\in(0,1)$ and $r_0=dist(0,\partial\Omega)$. Then\\
$(i)$ if $\max\{1+\frac{2\alpha}N,\frac{2\alpha}{N-2\alpha}\}<p<p^*_{\alpha}$,
there exist $R_0\in(0,r_0)$ and $c_{36}>0$ such that
\begin{equation}\label{3.1}
u_{\infty}(x)\ge c_{36}|x|^{-\frac{2\alpha}{p-1}}\qquad \forall x\in B_{R_0}\setminus\{0\},
\end{equation}
$(ii)$ if $\frac{2\alpha}{N-2\alpha}>1+\frac{2\alpha}N$ and $p=\frac{2\alpha}{N-2\alpha}$,
there exist $R_0\in(0,r_0)$ and $c_{37}>0$ such that
\begin{equation}\label{ying 3}
u_{\infty}(x)\ge \frac{c_{37}}{(1+|\log(|x|)|)^{\frac1{p-1}}}|x|^{-\frac{p(N-2\alpha)}{p-1}},\quad\forall  x\in B_{R_0}\setminus\{0\},
\end{equation}
$(iii)$ if $\frac{2\alpha}{N-2\alpha}>1+\frac{2\alpha}N$ and $p\in(1+\frac{2\alpha}N,\frac{2\alpha}{N-2\alpha})$,
there exist $R_0\in(0,r_0)$ and $c_{38}>0$ such that
\begin{equation}\label{3.2}
u_{\infty}(x)\ge c_{38}|x|^{-\frac{p(N-2\alpha)}{p-1}}\qquad\forall  x\in B_{R_0}\setminus\{0\}.
\end{equation}
\end{proposition}
{\bf Proof.} $(i)$
Using (\ref{1.3}) and Lemma \ref{lemma 3}$(i)$ with $\max\{1+\frac{2\alpha}N,\frac{2\alpha}{N-2\alpha}\}<p<p^*_{\alpha}$, we see that there exist $\rho_0\in(0,r_0)$ and $c_{39},c_{40}>0$ such that
\begin{equation}\label{3.2}
u_k(x)\ge c_{39}k|x|^{-N+2\alpha}-c_{40}k^p|x|^{-(N-2\alpha)p+2\alpha}\qquad \forall x\in B_{\rho_0}\setminus\{0\}.
\end{equation}
Set
\begin{equation}\label{3.3}
 \rho_k=(2^{(N-2\alpha)p-2\alpha-1}\frac{c_{40}}{c_{39}}k^{p-1})^{\frac1{(N-2\alpha)(p-1)-2\alpha}}.
\end{equation}
Since $(N-2\alpha)(p-1)-2\alpha<0$, there holds $\lim_{k\to\infty}\rho_k=0$. Let $k_0>0$ such that
$\rho_{k_0}\le \rho_0$,
then for $x\in B_{\rho_k}\setminus B_{\frac{\rho_k}2}$, we have
\begin{eqnarray*}
 c_{40}k^p|x|^{-(N-2\alpha)p+2\alpha}&\le& c_{40}k^p(\frac{\rho_k}2)^{-(N-2\alpha)p+2\alpha} \\
   &=& \frac{c_{39}}2k \rho_k^{-N+2\alpha} \\
   & \le& \frac{c_{39}}2k |x|^{-N+2\alpha}
\end{eqnarray*}
and
\begin{eqnarray*}
 k =(2^{(N-2\alpha)p-2\alpha-1}\frac{c_{40}}{c_{39}})^{-\frac1{p-1}}\rho_k^{N-2\alpha- \frac{2\alpha}{p-1}} \ge c_{41}|x|^{N-2\alpha- \frac{2\alpha}{p-1}},
\end{eqnarray*}
where $c_{41}=(2^{(N-2\alpha)p-2\alpha-1}\frac{c_{40}}{c_{39}})^{-\frac1{p-1}}2^{(N-2\alpha)(p-1)-2\alpha-1}$.
Combining with (\ref{3.1}), we obtain
\begin{eqnarray}
 u_k(x) &=& c_{39}k|x|^{-N+2\alpha}-c_{40}k^p|x|^{-(N-2\alpha)p+2\alpha} \nonumber\\
   &\ge&  \frac{c_{39}}2k|x|^{-N+2\alpha}\nonumber\\
   &\ge& c_{42}|x|^{-\frac{2\alpha}{p-1}},\label{3.7}
\end{eqnarray}
for $x\in B_{\rho_k}\setminus B_{\frac{\rho_k}2}$, where $c_{42}=\frac12c_{39}c_{41}$ is independent of $k$.
By (\ref{3.3}), we can choose a sequence $\{k_n\}\subset [1,+\infty)$ such that
$$\rho_{k_{n+1}}\ge \frac12\rho_{k_n},$$
For any $x\in B_{\rho_{k_0}}\setminus\{0\}$, there exists $k_n$ such that
$x\in B_{\rho_{k_n}}\setminus B_{\frac{\rho_{k_n}}2}$, then, by (\ref{3.7}),
$$u_{k_n}(x)\ge c_{42}|x|^{-\frac{2\alpha}{p-1}}.$$
Together with  $u_{k+1}>u_k$, we derive
$$u_{\infty}(x)\ge c_{42}|x|^{-\frac{2\alpha}{p-1}},\quad x\in B_{\rho_{k_0}}\setminus\{0\}.$$

\medskip

$(ii)$
By (\ref{1.3}) and Lemma \ref{lemma 3}-$(ii)$ with $p=\frac{2\alpha}{N-2\alpha}$, there exist $\rho_0\in(0,r_0)$ and $c_{43},c_{44}>0$ such that
\begin{equation}\label{3.2.2}
u_k(x)\ge c_{43}k|x|^{-N+2\alpha}-c_{44}k^p|\log(|x|)|,\quad x\in B_{\rho_0}\setminus\{0\}.
\end{equation}
Let $\{\rho_k\}$  be a sequence of real numbers with value  in $(0,1)$ and such that
\begin{equation}\label{3.2.3}
 c_{44}k^{p-1}|\log(\frac{\rho_k}2)|=\frac{c_{43}}2\rho_k^{-N+2\alpha}.
\end{equation}
Then $\lim_{k\to\infty}\rho_k=0$
and there exists $k_0>0$ such that
$\rho_{k_0}\le \rho_0$.
Thus, for any $x\in B_{\rho_k}\setminus B_{\frac{\rho_k}2}$ and $k\ge k_0$,
$$
c_{43}k^p |\log(|x|)|\le c_{44}k^p |\log(\frac{\rho_k}2)|=\frac{c_{43}}2k \rho_k^{-N+2\alpha}   \le \frac{c_{43}}2k |x|^{-N+2\alpha}.
$$
Therefore, assuming always $x\in B_{\rho_k}\setminus B_{\frac{\rho_k}2}$, we derive from (\ref{3.2.3}) that
$$
 k=(\frac{c_{44}}{2c_{43}})^{-\frac1{p-1}}(\frac{\rho_k^{-N+2\alpha}}{1+|\log(\rho_k)|})^{\frac{1}{p-1}}
 \ge c_{45}\frac{|x|^{- \frac{N-2\alpha}{p-1}}}{(1+|\log(|x|)|)^{\frac{1}{p-1}}},
$$
where $c_{45}=2^{ -\frac{N-2\alpha}{p-1}}(\frac{c_{44}}{2c_{43}})^{-\frac1{p-1}}$. Consequently
\begin{eqnarray}
 u_k(x) &\ge& c_{43}k|x|^{-N+2\alpha}-c_{44}k^p|\log(|x|)|\nonumber
   \\&\ge&  \frac{c_{43}}2k|x|^{-N+2\alpha}
   \ge c_{46}\frac{|x|^{-\frac{p(N-2\alpha)}{p-1}}}{(1+|\log(|x|)|)^{\frac{1}{p-1}}},\label{3.2.7}
\end{eqnarray}
where $c_{46}=\frac12 c_{43}c_{45}$ is independent of $k$.\smallskip

By (\ref{3.2.3}), we can choose a sequence $k_n\in [1,+\infty)$ such that
$$\rho_{k_{n+1}}\ge \frac12\rho_{k_n},$$
Then for any $x\in B_{\rho_{k_0}}\setminus\{0\}$, there exists $k_n$ such that
$x\in B_{\rho_{k_n}}\setminus B_{\frac{\rho_{k_n}}2}$.  By (\ref{3.2.7}) there holds
$$u_{k_n}(x)\ge c_{46}\frac{|x|^{-\frac{p(N-2\alpha)}{p-1}}}{(1+|\log(|x|)|)^{\frac{1}{p-1}}}.$$
Together with  $u_{k+1}>u_k$, we infer
$$u_{\infty}(x)\ge c_{46}\frac{|x|^{-\frac{p(N-2\alpha)}{p-1}}}{(1+|\log(|x|)|)^{\frac{1}{p-1}}}\qquad \forall x\in B_{\rho_{k_0}}\setminus\{0\}.$$

\medskip
$(iii)$ By (\ref{1.3}) and Lemma \ref{lemma 3}-$(iii)$ with $p\in(1+\frac{2\alpha}N,\frac{2\alpha}{N-2\alpha})$, there exist $\rho_0\in(0,r_0)$ and $c_{47},c_{48}>0$ such that
\begin{equation}\label{3.1.2}
u_k(x)\ge c_{47}k|x|^{-N+2\alpha}-c_{48}k^p\qquad\forall x\in B_{\rho_0}\setminus\{0\}.
\end{equation}
Put
\begin{equation}\label{3.1.3}
 \rho_k=(\frac{c_{48}}{2c_{47}}k^{p-1})^{-\frac1{N-2\alpha}},
\end{equation}
then $\lim_{k\to\infty}\rho_k=0$
and there exists $k_0>0$ such that
$\rho_{k_0}\le \rho_0$.
Therefore, if $x\in B_{\rho_k}\setminus B_{\frac{\rho_k}2}$ and $k\ge k_0$, there holds
$$
 c_{48}k^p   =\frac{c_{47}}2k \rho_k^{-N+2\alpha}   \le \frac{c_{47}}2k |x|^{-N+2\alpha},
$$
which yields
$$
 k=(\frac{c_{48}}{2c_{47}})^{-\frac1{p-1}}\rho_k^{ -\frac{N-2\alpha}{p-1}}
 \ge c_{49}|x|^{- \frac{N-2\alpha}{p-1}},
$$
by (\ref{3.1.3}), where $c_{49}=2^{ -\frac{N-2\alpha}{p-1}}(\frac{c_{48}}{2c_{47}})^{-\frac1{p-1}}$.
Consequently,  
\begin{eqnarray}
 u_k(x) \ge c_{47}k|x|^{-N+2\alpha}-c_{48}k^p
   &\ge&  \frac{c_{47}}2k|x|^{-N+2\alpha}\nonumber\\
   &\ge& c_{50}|x|^{-\frac{p}{p-1}(N-2\alpha)},\label{3.1.7}
\end{eqnarray}
where $c_{50}=\frac12 c_{47}c_{49}$ is independent of $k$.

By (\ref{3.1.3}), we can choose a sequence $k_n\in [1,+\infty)$ such that
$$\rho_{k_{n+1}}\ge \frac12\rho_{k_n},$$
Then for any $x\in B_{\rho_{k_0}}\setminus\{0\}$, there exists $k_n$ such that
$x\in B_{\rho_{k_n}}\setminus B_{\frac{\rho_{k_n}}2}$ and then by (\ref{3.1.7}),
$$u_{k_n}(x)\ge c_{50}|x|^{-\frac{p(N-2\alpha)}{p-1}}.$$
Together with  $u_{k+1}>u_k$, we have
$$u_{\infty}(x)\ge c_{50}|x|^{-\frac{p(N-2\alpha)}{p-1}}\qquad \forall x\in B_{\rho_{k_0}}\setminus\{0\},$$
which ends the proof.\hfill$\Box$
\begin{lemma}\label{upper bound}
Let $p\in(1+\frac{2\alpha}N,p^*_\alpha)$ and $u_s$ be a strongly singular solution of (\ref{eq1.2}) satisfying (\ref{eq1.01}). Then
\begin{equation}\label{4.1}
 u_{\infty}\le u_s\quad{\rm in}\quad \Omega\setminus\{0\},
\end{equation}
where $u_\infty$ is defined by (\ref{definition infty}).
\end{lemma}
{\bf Proof.} By (\ref{eq1.01}) and (\ref{1.3}),
it follows 
$$\lim_{x\to0} u_s|x|^{\frac{2\alpha}{p-1}}=c_0\quad
{\rm and} \quad \lim_{x\to0} u_k|x|^{N-2\alpha}=c_k,$$
which implies that there exists $r_1>0$ such that
$$u_k<u_s\quad{\rm in}\quad B_{r_1}\setminus\{0\}.$$
Since by Theorem \ref{pr 2.1}, $u_k$ satisfies
$$
(-\Delta)^\alpha u_k+u_k^p=0\quad{\rm in}\quad\Omega\setminus B_{r_1}(0),
$$
 so does $u_s$. By Theorem \ref{teo CP} there holds
 $u_k\le u_s$ in $\Omega\setminus\{0\}$. Jointly with (\ref{definition infty}), it implies
$$u_{\infty}\le u_s\quad{\rm in}\quad \Omega\setminus\{0\}.$$
 \hfill$\Box$

\medskip
\noindent{\bf Proof of Theorem \ref{teo 1} $(ii)$ and Theorem \ref{teo 1.1} $(iv)$.}
By Lemma \ref{upper bound} and Theorem \ref{teo 2}, we  obtain that
$u_{\infty}$ is a classical solution of (\ref{eq1.2}). Moreover, by Proposition \ref{pr 3.2} part $(i)$ and Lemma \ref{upper bound},
we have
$$\frac1{c_{51}}|x|^{-\frac{2\alpha}{p-1}}\le u_{\infty}(x)\le c_{51}|x|^{-\frac{2\alpha}{p-1}},$$
for some $c_{51}>1$.
Then $u_{\infty}=u_s$ in $\R^N\setminus\{0\}$ since $u_s$ is unique in the class of solutions satisfying (\ref{strong singular unique}). 
\hfill$\Box$
\subsection{Proof of Theorem \ref{teo 1.1} $(ii)$ and $(iii)$ }
 By Lemma \ref{upper bound} and Theorem \ref{teo 2},
$u_{\infty}$ is a classical solution of (\ref{eq1.2}) and it satisfies
$$u_\infty\le u_s\quad {\rm in}\quad \Omega\setminus\{0\}.$$
Therefore (\ref{1.1.1}) and (\ref{1.2.1}) follow by Proposition \ref{pr 3.2} part $(ii)$ and $(iii)$, respectively.

\hfill$\Box$


\end{document}